\documentclass[12pt]{article}
\usepackage{amssymb,amsmath,makeidx, amsthm}
\usepackage{enumerate}
\usepackage{psfig}
\usepackage{color}
\usepackage[all]{xy}
\begin{document}
\newtheorem{theorem}{Theorem}[subsection]
\newtheorem{pusto}[theorem]{}
\numberwithin{subsection}{section}
\newtheorem{remark}[theorem]{Remark}
\newcommand{\po}{{\perp_{\omega}}}
\newcommand{\mathto}{\mathop{\longrightarrow}\limits}
\def\bt{\begin{theorem}}
\def\et{\end{theorem}}
\def\x{\index}
\def\bp{\begin{pusto}}
\def\ep{\end{pusto}}
\def\bi{\begin{itemize}}
\def\ei{\end{itemize}}
\def\btl{$\blacktriangleleft$}
\def\btr{$\blacktriangleright$}
\def\mn{\medskip\noindent}
\def\n{\noindent}
\def\p{\partial}
\def\e{\varepsilon}
\def\a{\alpha}
\def\d{\delta}
\def\s{\sigma}
\def\ol{\overline}
\def\wt{\widetilde}
\def\wh{\widehat}
\def\Cal{\cal}
\def\cb{\color{blue}}
\def\Imm{\mathrm{Imm}}
\def\Int{\mathrm{Int}}
\def\Ker{\mathrm{Ker}}
\def\Id{\mathrm{Id}} 
\def\id{\mathrm{id}}
\def\su{\mathrm{supp}}
\def\rank{\mathrm{rank}}
\def\dim{\mathrm{dim}}
\def\codim{\mathrm{codim}}
\def\const{\mathrm{const}}
\def\Fix{\mathrm{Fix}}
\def\Span{\mathrm{Span}}
\def\span{\mathrm{span}}
\def\inv{\mathrm{inv}}
\def\Hom{\mathrm{Hom}}
\def\Mono{\mathrm{Mono}}
\def\Diff{\mathrm{Diff}}
\def\Emb{\mathrm{Emb}}
\def\emb{\mathrm{emb}}
\def\imm{\mathrm{imm}}
\def\sub{\mathrm{sub}}
\def\Hol{\mathrm{Hol}}
\def\hol{\mathrm{hol}}
\def\Sec{\mathrm{Sec}}
\def\Sol{\mathrm{Sol}}
\def\Supp{\mathrm{Supp}}
\def\bs{\mathrm{bs}}
\def\Vect{\mathrm{Vect}}
\def\Vert{\mathrm{Vert}}
\def\Clo{\mathrm{Clo}}
\def\clo{\mathrm{clo}}
\def\Exa{\mathrm{Exa}}
\def\Gr{\mathrm{Gr}}
\def\CR{\mathrm{CR}}
\def\Distr{\mathrm{Distr}}
\def\dist{\mathrm{dist}}
\def\Ham{\mathrm{Ham}}
\def\mers{\mathrm{mers}}
\def\Iso{\mathrm{Iso}}
\def\iso{\mathrm{iso}}
\def\tang{\mathrm{tang}}
\def\trans{\mathrm{trans}}
\def\ttr{\mathrm{trans-tang}}
\def\ot{\mathrm{ot}}
\def\loc{\mathrm{loc}}
\def\Symp{\mathrm{Symp}}
\def\symp{\mathrm{symp}}
\def\Lag{\mathrm{Lag}}
\def\isosymp{\mathrm{isosymp}}
\def\subisotr{\mathrm{sub-isotr}}
\def\Cont{\mathrm{Cont}}
\def\cont{\mathrm{cont}}
\def\Leg{\mathrm{Leg}}
\def\isocont{\mathrm{isocont}}
\def\isot{\mathrm{isot}}
\def\isotr{\mathrm{isotr}}
\def\coisot{\mathrm{coisot}}
\def\coreal{\mathrm{coreal}}
\def\real{\mathrm{real}}
\def\comp{\mathrm{comp}}
\def\old{\mathrm{old}}
\def\new{\mathrm{new}}
\def\top{\mathrm{top}}
\def\rG{\mathrm{G}}
\def\CS{\mathrm{CS}}
\def\ND{\mathrm{ND}}
\def\CND{\mathrm{CND}}
\def\fA{\mathfrak {A}}
\def\fa{\mathfrak {a}}
\def\A{\mathcal {A}}
\def\E{\mathcal {E}}
\def\F{\mathcal {F}}
\def\G{\mathcal {G}}
\def\H{\mathcal {H}}
\def\I{\bf {I}}
\def\J{\mathcal {J}}
\def\N{\mathcal {N}}
\def\R{\mathcal {R}}
\def\P{\mathcal {P}}
\def\T{\mathcal {T}}
\def\L{\mathcal {L}}
\def\SS{\mathcal {S}}
\def\U{\mathcal {U}}
\def\D{\mathcal {D}}
\def\V{\mathcal {V}}
\def\W{\mathcal {W}}
\def\Z{\mathcal {Z}}
\def\C{\mathcal {C}}
\def\O{\mathcal {O}}
\def\BH{{\bold H}}
\def\Bf{ f} 
\def\Bdelta{{\delta}}
\def\BF{{\bold F}}
\def\BS{{\bold S}}
\def\bbZ{{\mathbb Z}}
\def\bbR{{\mathbb R}}
\def\bbS{{\mathbb S}}
\def\bbC{{\mathbb C}}
\def\bbQ{{\mathbb Q}}
\def\IS{\R_{\sympl-iso}}
\def\IC{\R_{\cont-iso}}
\def\Sp{\mathrm{Sp}}
\def\GL{\mathrm{GL}}
\def\SO{\mathrm{SO}}
\def\U{\mathrm{U}}
\def\X{X^{(r)}}
\def\oC{\overline{C}}
\newcommand{\eps}{\epsilon}
\def\Op{{\mathcal O}{\it p}\,}
\def\hook{\lrcorner\,}
\newcommand{\fM}{{\frak M}}
\newcommand{\fm}{{\frak m}}
\newcommand{\fI}{{\frak I}}
\centerline{\LARGE Topology of spaces of  $S$-immersions}
\medskip
 \centerline
{\small To Oleg Viro on his 60th birthday}
\bigskip
\bigskip
\qquad   \qquad \quad 
Y. M. Eliashberg
\footnote{Partially supported by NSF 
grants DMS-0707103 and DMS 0244663}
\qquad \qquad \quad N. M. Mishachev

\quad \qquad \qquad Stanford University \qquad Lipetsk Technical
University


\begin{abstract}
\noindent
We use the wrinkling theorem proven in 
\cite{[EM97]} to  fully describe the homotopy type of the space of 
$S$-immersions, i.e. equidimensional folded maps with prescribed folds. 
\end{abstract}

\tableofcontents
 
\section{Equidimensional folded maps}
\label{s:main}

Let $V$ and $W$ be  two manifolds of dimension  $q$.
The manifold $V$ will always be  assumed   {\it closed} and {\it connected.}
A map $f:V\to W$ is called {\it folded},
if it has only fold type singularities.
We will discuss in this paper an $h$-principle for folded mappings
$f:V\to W$ with a  {\it prescribed fold} $\Sigma^{10}(f)\subset V$.
Folded maps with the fold $S=\Sigma^{1,0}(f)$ are also called
{\it $S$-immersions}, see \cite{[El70]}. {\bf Throughout this paper we assume 
that $S\neq\varnothing$.}
We provide in this paper an essentially complete description of the homotopy 
type of the spaces of $S$-immersions.
More precisely, we prove that in most cases one has  an $h$-principle type 
result, while sometimes the  space of $S$-immersions may have a number of 
additional components of a different nature. The topology of these components is  also  
fully described.

{\bf Remarks.}
{\bf 1.} The results proven in this paper were formulated in  \cite{[El70],[El72]}, 
but the proof of the  injectivity part of the $h$-principle claim was never 
published before. 

{\bf 2.} The approach described in this paper to the problem of construction of 
mappings with prescribed singularities generalizes to the case of maps $V^n\to 
W^q$ for $n>q$.  However,  unlike the case $n=q$,  the results which can be 
proven using the  current techniques  are essentially equivalent  to the results 
of \cite{[El72]}.

{\bf 3.} The subject of this paper is the geometry   of singularities in the {\it 
source} manifold.
The geometry of singularities in the {\it image} is much
more subtle, see \cite{[Gr07]}.


\subsection{$S$-immersions}
\label{ss:s-immersions}

We refer the reader to Section \ref{s:wrinkles} below for the definition
and basic properties of fold, wrinkles, and for the formulation of
the Wrinkling Theorem from \cite{[EM97]}.

Let $f:V\to W$ be an $S$-immersion, i.e.  a map  with only
fold type singularity $S=\Sigma^{1,0}(f)$.
The fold $S\subset V$ has a neighborhood $U$ which admits an involution
$\alpha_{\,\loc}:U\to U$ such that  
$f\circ\alpha_\loc=f$. 
In particular, if $S$ divides $V$  into two submanifolds $V_\pm$ 
with the common boundary $\p V_\pm=S$ 
then an $S$-immersion is just a pair of immersions $f_\pm:V_\pm\to W$ such
that $f_+=f_-\circ\alpha_{\,\loc}$ near $S$.
Denote by $T_SV$ an  $n$-dimensional tangent bundle over $V$
which is obtained from $TV$ by re-gluing $TV$ along $S$
with $d\alpha_\loc$. For example, if $V=S^q$ and $S$ is the equator
$S^{q-1}\subset S^q$, then $T_SV=S^q\times \bbR^q$.
We will call $T_SV$ the {\it tangent bundle of $\,V$
folded along $S$}. The differential $df: TV\to TW$ of any $S$-immersion
$f:V\to W$ has a canonical (bijective) regularization 
$d_S f:T_SV\to TW$, the {\it folded  differential} of $f$.

\mn 
Let us denote by $\fM(V,W,S)$ the space of $S$-immersions $V\to W$.
We also consider the space $\fm(V,W,S)$, a formal analog of $\fM(V,W,S)$
which consists of bijective homomorphisms $T_SV\to TW$.
The folded differential induces a natural inclusion $d:\fM(V,W,S)\to \fm (V,W, 
S)$.
Our goal is to study the homotopical properties of the map $d$. We prove that
in most cases the map $d$ is a (weak) homotopy equivalence. However, there are 
some 
exceptional cases when the map $d$ is a homotopy equivalence on some of
the components of $ \fM(V,W,S)$, while the structure of remaining components
can also be completely understood.
  
\subsection{Taut-soft dichotomy for  $S$-immersions} 
\label{ss:taut}

A map $f\in\fM(V,W,S)$ is called {\it taut}
if there exists an involution $\alpha:V\to V$ such that
${\rm Fix}\, \alpha= S$ and $f\circ\alpha=f$.
We denote by ${\frak M}_{taut}(V,W,S)$ the subspace of $\fM(V,W,S)$
which consists of taut maps. Non-taut maps are called {\it soft}, and  we denote 
${\frak M}_{soft}(V,W,S):= \fM(V,W,S)\setminus\fM_{taut}(V,W,S)$.
A map $f\in{\frak M}_{taut}(V,W,S)$ uniquely determines the corresponding
involution $\alpha$, and thus 
\bp {\bf (Topological structure of the space of taut $S$-immersions)}
The space $ {\frak M}_{taut}(V,W,S)$ is the space of pairs $(\alpha,h)$
where  $\alpha:V\to V$ is an involution with ${\rm Fix}\, \alpha= S$
and $h$ is an immersion of the  quotient manifold $V/\alpha$ with
the boundary $S$ to $W$.
\ep
Of course, for most pairs $(V,S)$ 
the space $\fI(V,S)$ of such involutions  is empty, and hence
in these cases the space ${\frak M}_{taut}(V,W,S)$ is empty as well.  

\mn
The topology of the space ${\frak M}_{taut}(V,W,S)$ is especially simple if 
$S$ divides $V$ into two submanifolds $V_\pm$ with the common boundary
$\p V_\pm=S$, e.g when both manifolds $V$ and $S$ are orientable. Clearly, 

\bp
\label{t:taut-fibration}
{\bf (Topological structure of the space of taut $S$-immersions: the orientable case)}
If $S$ divides $V$  into two submanifolds $V_\pm$ with the common
boundary $\p V_\pm=S$, such that there exists a diffeomorphism
$V_+\to V_-$ fixed along the boundary, then the space
${\frak M}_{taut}(V,W,S)$ is homeomorphic
to the product 
$$\Diff_S(V_+)\times\Imm(V_+,W)\,,$$
where $\Diff_S(V_+)$
is the group of diffeomorphisms $V_+\to V_+$ fixed at the boundary
together with their $\infty$-jet, and $\Imm(V_+,W)$
is the space of immersions $V_+\to  W$.
\ep

\n
Note that according to Hirsch's theorem, \cite{[Hi]},  the space  
$\Imm(V_+,W)$ is homotopy equivalent to the space $Iso(V,W)$ of fiberwise
isomorphic bundle maps $TV_+\to TW$.
For instance, when $V=S^q, W=\bbR^q$ and $S$ is the equator
$S^{q-1}\subset S^q$ then we get 
$${\frak M}_{taut}(S^q,\bbR^q,S^{q-1})\mathop{\simeq}\limits^{h.e.}
\Diff_{\p D^q} D^q\times O(q).$$
 In particular, when $q=2$ the space  of taut  $S^1$-immersions
$S^2\to\bbR^2$ which preserve  orientation on $S^2_+$ is homotopy equivalent to 
$S^1$. 

\bp
\label{p:taut}
{\bf(Subspaces of taut and soft $S$-immersions are open-closed)}
Let $S\subset V$ be a closed $(q-1)$-dimensional submanifold of $\,V$.
Then the subspaces ${\frak M}_{taut}(V,W,S)\subset\fM(V,W,S)$ and
${\frak M}_{soft}(V,W,S)$  are open
and closed, i.e.  they consist of whole connected components of $\fM(V,W,S)$
\ep
\begin{proof}
Given any $f\in\fM(V,W,S)$ there exists an $\e=\e(f)>0$ such that  the local 
involution 
$\alpha^f_\loc:\Op S\to\Op S$ with  $f\circ\alpha^f_\loc=f$ is defined on
an $\e$-tubular neighborhood   $U_\e\supset S$.  
If $f$ is taut then  $\alpha^f_\loc$ extends as a global involution  
$\alpha^f:V\to V$ such that $f\circ\alpha^f=f$.   Hence, for any $x\in 
V\setminus U_{\e/2}$ the distance $d(x,\alpha^f(x))\geq  \e$.  This implies that  there 
exists a $\delta(\e)>0$ such that  if $||f'-f||_{C^2}<\delta(\e)$   then 
$\e(f')>\frac{\e}2$.  Therefore the local involution  $\alpha_\loc^{f'}$ is 
defined on $U_{\e(f)/2}$. We extend   it  to $V\setminus U_{\e(f)/2}$ as  a 
global involution $\alpha^{f'}$ by defining  $\alpha^{f'}(x)$ as the unique 
point from $(f')^{-1}(x)$ whose distance from $\alpha^f(x)$ is 
$<\frac{\e(f)}2$. Hence, $\fM_{taut}(V,W,S)$ is open.
 On the other hand, using  the same argument together with the implicit function 
theorem we conclude  that   given a sequence $f_n\in \fM_{taut}(V,W,S)$,
such that  $f_n\mathop{\to}\limits^{C^2} f$ then
$\alpha^{f_n}\mathop{\to}\limits^{C^1}\alpha^f$, and hence $f\in \fM_{taut}(V,W,S)$
and $\fM_{taut}(V,W,S)$ is closed.
\end{proof}


\medskip
The following theorem completes the description of the homotopy type 
of the spaces of $S$-immersions.

\bp
\label{t:soft}
{\bf (Homotopical structure of the space of soft $S$-immersions)}
Let $S\subset V$ be a closed non-empty $(q-1)$-dimensional submanifold of $V$. 
Suppose that  $q\geq3$, or $q=2$ but $W$ is open. Then the inclusion 
$$d: {\frak  M}_{soft}(V,W,S)\to {\frak m}(V,W,S)$$
is a (weak) homotopy equivalence. In particular, when the space $\fI(V,S)$
is empty (and hence ${\frak M}_{taut}(V,W,S)$ is empty as well) then
$d: {\frak M} (V,W,S)\to \frak m(V,W,S)$ is a (weak) homotopy equivalence.
\ep

\begin{figure}
\centerline{\psfig{figure=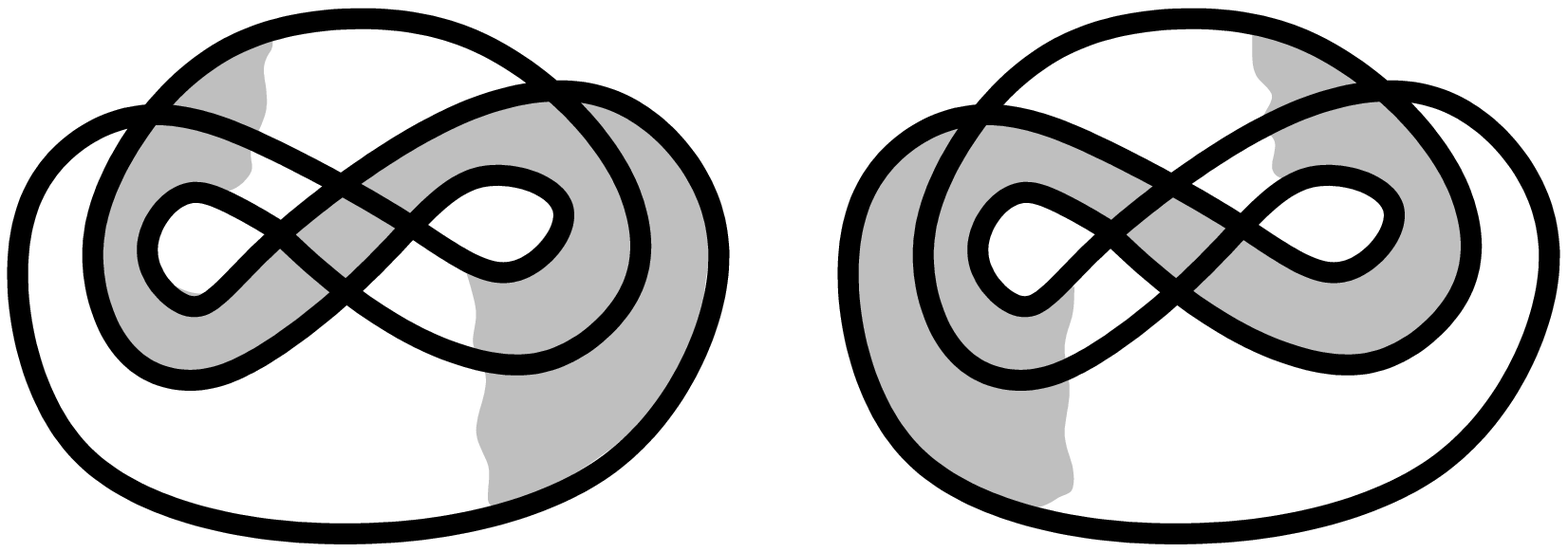,height=30mm}}
\caption{${f\in \frak M}_{soft} (S^2,\bbR^2,S^1)$ as a pair $f_1,f_2: D^2\to 
\bbR^2$}
\label{b02}
\end{figure}

For instance, when $V=S^q, W=\bbR^q$ and $S$ is the equator
$S^{q-1}\subset S^q$ then we get 
$${\frak M}_{soft}(S^q,\bbR^q,S^{q-1})\mathop{\simeq}\limits^{h.e.}
\Omega_q(O(q)),$$
where  $\Omega_q(O(q))$ is the (free) $q$-loops space of the orthogonal group
$O(q)$.  In particular,  using also   \ref{t:taut-fibration} we conclude 
that   the space  of {\it all} $S^1$-immersions
$S^2\to\bbR^2$ which preserve orientation on $S^2_+$ consists of two components  
homotopy equivalent to $S^1$.  The standard projection represents the taut component while  
an example of a soft $S^1$-immersion $f:S^2\to \bbR^2$
(as a pair of maps $f_1,f_2:D^2\to\bbR^2$, $f_1|_{S^1}= f_2|_{S^1}$)
is presented on Fig.\ref{b02}.   

\mn
{\bf Remarks.}
{\bf 1.} Theorem  \ref{t:soft} 
was formulated in \cite{[El72]}, but only the epimorphism part of  the statement
was proven there.  

{\bf 2.} Theorem  \ref{t:soft}  trivially holds for $q=1$ and any $W$.

{\bf 3.} In the case $q=2$ and $W$ is a closed surface we can only prove that 
the map $d$ induces an epimorphism on homotopy groups.

\section{Zigzags}
\label{s:zigzags}

\subsection{Zigzags and soft $S$-immersions}
\label{ss:zigzags}
A {\it model zigzag} is any smooth function $\Z:[a,b]\to\bbR$ such that
\begin{itemize}
\item $\Z$ is increasing near $a$ and $b$;
\item $\Z$ has exactly two interior non-degenerate critical points 
$m,M\in(a,b)$, $m<M$;
\item $\Z(b)>\Z(a)$.
\end{itemize} 
If, in addition, $\Z(b)>\Z(M)$ and $\Z(a)<\Z(m)$ then the model zigzag is called   
{\it long}. 

{\bf Example.} The function $\Z(z)=z^3-3z$ is a model zigzag on  any interval
$[a,b]$ such that $a<-1,b>1$ and $\Z(b)>\Z(a)$. If $a<-2, b>2$ then $\Z$ is 
long.

An embedding $h:[a,b]\to V$ is called a {\it zigzag}  (resp. {\it long zigzag}) 
of a map $f\in\fM(V,W,S)$ if 
\begin{itemize}
\item it is transversal to $S$, and
\item the composition $f\circ h$ can be presented as
$g\circ \Z$, where $\Z:[a,b]\to\bbR$ is a model zigzag (resp. long model zigzag) 
and $g:[A,B]\to W$ is an immersion defined on an interval $[A,B]$ such that 
$(A,B)\supset\Z([a,b])$.
\end{itemize}
\begin{figure}
\centerline{\psfig{figure=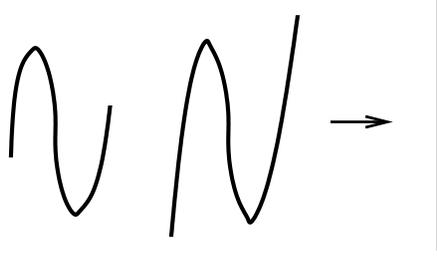,height=35mm}}
\caption{Zigzag and long zigzag.}
\label{em41}
\end{figure}
The image $h([a,b])$ of the embedding $h$ will also be called a (long) 
zigzag. The points $h(a),h(b)\in V$ are called  the {\it end points} of the zigzag 
$h$. Note that
\bp
\label{l:ext}
{\bf (Extension of long zigzags)}
Let $h:[a,b]\to V$ be a {\bf long} zigzag. Then any extension $h':[a',b']\to V$
of the embedding $h$ which does not have additional intersection points with $S$ 
is a long zigzag.
\ep
The implicit function theorem implies:
\bp
\label{l:flex}
{\bf (Local flexibility of zigzags)}
Let  $f_t\in\fM(V,W,S)$ defined for $t\in\Op0\subset\bbR$. Suppose that $f_0$ 
admits a zigzag
$h_0:[a,b]\to V$ such that the composition $f_0\circ h_0$  can be   factored as 
$g_0\circ\Z$ for an immersion $g_0:[A,B]\to V$. Suppose that $g_0$ is included 
into a family of immersions   $g_t:[A,B]\to W$, $t\in\Op 0$. Then there exists a 
family of zigzags $h_t:[a,b]\to V$ defined for  $t\in\Op 0$  such that
$f_t\circ h_t=g_t\circ\Z$.
\ep
 
\bp
\label{l:soft-zigzag}
{\bf (Softness criterion)}
A map $f\in\fM(V,W,S)$ is soft if and only if it admits a zigzag.
\ep
\begin{proof}
Existence of any zigzag is incompatible with the  existence of an involution, and 
hence a map admitting a zigzag is soft.
On the other hand, if $f\in\fM(V,W,S)$ does not admit a zigzag we can define   
  an involution $\alpha=\alpha^f:V\to V\,$,  $f\circ \alpha=f$, as follows.
 Given  $v\in V$ take any arc $C\subset V$ connecting $v$ in $V\setminus S$ with 
a point 
$s\in S$. 
Then the absence of zigzags guarantees that  $f^{-1}(C)$ contains a unique 
candidate $v'$  for $\alpha(v)$. Similarly, if for another path we had another 
candidate $v''$ this would create a zigzag. Hence, the involution   $\alpha 
:V\to V\,$ with  $f\circ \alpha=f$ is correctly defined, and therefore  $f$ is 
taut.
\end{proof}

\subsection{Zigzags adjacent to a chamber}\label{s:ample-z}

The components of $V\setminus S$ are called {\it chambers}.
We say that a zigzag $h$ is {\it adjacent to a chamber} $C$ if this 
chamber contains one of the end points of the zigzag.  

\mn 
A family of maps $f_s\in\fM_{soft}(V,W,S)$, $s\in K$, parameterized by a connected compact 
set
$K$ can be viewed as a fibered map $\wt f:K\times V\to K\times W$. We will call
such maps {\it fibered (over $K$)} $S$-immersions. A fibered $S$-immersion is 
called
{\it soft} if it is fiberwise soft.  According to \ref{p:taut} $\wt f$ is soft 
if it
is soft over a point $s\in K$. We denote the space of soft 
$S$-immersions fibered over $K$ by $\fM_{soft}^K(V,W, S)$.
A family $Z_s$ of zigzags for $f_s$, $s\in K$, will be referred to 
as a {\it fibered} over $K$ zigzag $Z$. A fibered zigzag is called {\it special} 
if  the projections
$f_s (Z^s)$ are independent of $s\in K$.
We say that a fibered soft map $\wt f:K\times V\to K\times W$ admits
a {\it set of fibered ({\rm resp.}  special fibered) zigzags subordinated to a 
covering}   $K=\mathop{\bigcup}\limits_1^N U_j$ 
if there exist  fibered (resp. special fibered) {\it disjoint}  zigzags
$\wt Z_{1}, \dots, \wt Z_{N}$ for the fibered maps
$\wt f_j=\wt f|_{U_j\times V}:U_j\times V\to U_j\times W, j=1,\dots,N$.
\bp
\label{l:inscribed}
{\bf(Set of zigzags subordinated to an inscribed 
covering)}  Let  $\wt Z_{1}, \dots, \wt Z_{N}$ be a set of (special) fibered 
zigzags for a fibered map $\wt f$ subordinated to a covering  
$K=\mathop{\bigcup}\limits_1^N U_j$. Let  $K=\mathop{\bigcup}\limits_1^{N'} 
U'_j$  be another covering inscribed to the first one, i.e. for every $U'_i$, 
$i=1,\dots, N'$, there is $U_j$, $j=1,\dots, N$, such that
$U_i'\subset U_j$.
Then there exists a set of fibered  (special) zigzags $\wt Z_1',\dots,\wt 
Z'_{N'}$ subordinated to the covering $K=\mathop{\bigcup}\limits_1^{N'} U'_j$ 
such that for each   $s\in U_j'$ the zigzag ${Z'}_j^s $ is $C^\infty$-close to 
the zigzag $Z^s_j$  for some $j=1,\dots, N$.
\ep

\begin{proof} Lemma \ref{l:flex} implies  that a neighborhood of any (special)
fibered zigzag is foliated by (special) fibered zigzags, and hence  near any
(special) fibered zigzag $\wt Z$ one can always
find an arbitrarily many
disjoint copies of    $\wt Z$.
\end{proof}
We would like  to prove that a fibered soft map admits a set of  fibered zigzags
adjacent to any of its chambers $C$. Below we  explain two methods for proving 
this.
However, each of the methods requires some additional assumptions. The first one works
for any $W$, but only if $q\geq 3$, while the second works for  $q\geq 2$,
but only if the  target manifold is open. We do not know whether the statement
still holds if $W$ is a closed surface.
\footnote{The second approach also works  for $q=1$ and any $W$. 
The case $W=T^2$  can also  be treated by a slight modification of 
this method.}
\bp
\label{l:ample-parametric}
{\bf (Fibered zigzags adjacent to  a chamber)}
Let $\wt f$ be a fibered over $K$ soft $S$-immersion.  Suppose that  $q\geq3$, 
or $q=2$ but $W$ is open. 
Then, given any chamber $C\subset V\setminus S$ there is a homotopy
$\wt f_t\in\fM_{soft}^K (V,W,S)$, $t\in[0,1]$, such that $\wt f_0=\wt f$
and $\wt f_1$ admits a  set of special fibered zigzags adjacent to the chamber 
$C$.
\ep
\begin{figure}
\centerline{\psfig{figure=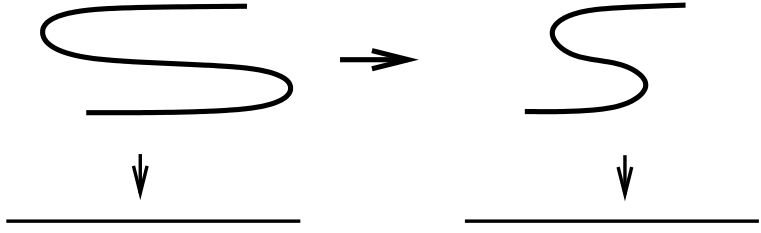,height=30mm}}
\caption{Making zigzag long}
\label{em11}
\end{figure}
{\sc Case $q\geq 3$.} We begin with two lemmas.
\bp
\label{l:long}
{\bf (Making zigzags long)}
Let $\wt f\in\fM^K (V,W,S)$. Suppose that  $q\geq3$.
Then there exists a homotopy  $\wt f_t\in\fM^K (V,W,S)$, $t\in[0,1]$, 
and a covering $K=\mathop{\bigcup}\limits_1^N U_j$ such that $\wt f_0=\wt f$
and $\wt f_1$ admits a set of special fibered long  zigzags subordinated
to a  covering  inscribed in the
covering $K=\mathop{\bigcup}\limits_1^N U_j$.
\ep

\begin{proof}
Choose a zigzag $h_s:[a,b]\to V$ for each $f_s$, $s\in K$,
in the family $\wt f$. Denote by $g_s$ the corresponding
immersion $[A,B]\to W$. According to Lemma \ref{l:flex}  
there is a neighborhood $U_s\ni s$ in $K$ such that for
each $s'\in U_s$ there exists a zigzag $h'_{s'}:[a,b]\to V$ such that
it factors through the same immersion $g_s$.
Due to compactness of $K$ we can choose a finite
subcovering $U_j=U_{s_j}$, $j=1,\dots, N$.
Lemma \ref{l:flex}  further implies that if we  $C^\infty$-small
perturb the immersions $g_j:=g_{s_j}$  then there still exist
for each $j=1,\dots N$, and  $s\in U_j$,  zigzags $Z_{s,j}\subset V$ 
which factor through them. Hence, if $q\geq 3$ the general position argument 
allows us to assume that the images  of the immersions $g_j$, 
and hence of the 
fibered zigzags $\wt Z_j$ do not intersect.
Now for each $j=1,\dots N$, there exists a deformation
$f_{s,t}$, $s\in U_j$, $t\in[0,1]$, of the family
$f_{s}$  supported  for each $s$ in an arbitrarily small
neighborhood of  $Z_s$ ,which makes the zigzags long, see Fig.\ref{em11}.
But the images of constructed zigzags do not intersect, and hence
if the neighborhoods of zigzags are chosen sufficiently small, then 
the above deformation can be done simultaneously over all $U_j$, $j=1,\dots, N$.
\end{proof}

\begin{figure}
\centerline{\psfig{figure=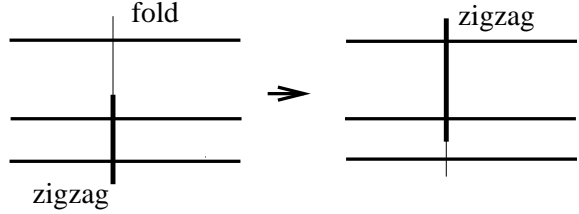,height=30mm}}
\caption{Penetration through a wall}
\label{em91}
\end{figure}

\begin{figure}
\centerline{\psfig{figure=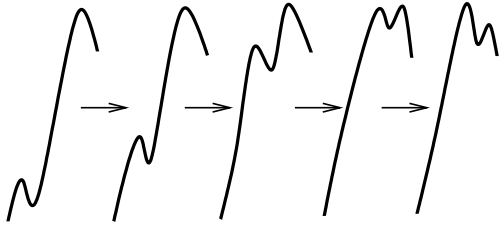,height=30mm}}
\caption{Homotopy $f_t$ on $h([a,c'])$}
\label{em24}
\end{figure}

\bp
\label{l:moving-long}
{\bf (Penetration through walls)}
Let  $C$ be one of the chambers for $f\in\fM(V,W,S)$.
Let $h:[a,c]\to V$ be an embedding such that
\begin{itemize}
\item there exists $b\in(a,c)$ such that $h|_{[a,b]}$ is a long zigzag for $f$;
\item $h$ is transversal to $S$ and $h(c)\in C$.
\end{itemize}
Then there exists a deformation $f_t\in\fM(V,W,S)$, $t\in[0,1]$,  with  $f_0=f$
which is supported in the neighborhood of the image $h([a,c])$ and such that
for some $b'\in(a,c)$ the embedding 
$h|_{[b',c]}$ is a zigzag for $f_1$ adjacent to $C$.
\ep
\begin{proof}
Suppose that the embedding $h|_{[b,c]}$ intersects the wall $S$ in a 
sequence of points $h(p_1),\dots h(p_k)$. We can push the original zigzag 
consequently through these points. Let $k=1$. By a small perturbation of $f$ 
near $h([a,c])$ we can make $h$ invariant with respect to the local 
involution $\alpha_{loc}$ on $\Op p$. Then there exists $c'\in (p,c)$
such that $f\circ h|_{[a,c']}$ can be factored as $[a,c']\to \bbR\to W$.
Then the deformation shown on  Fig.\ref{em91} and 
Fig.\ref{em24} allows  to move  the zigzag from $[a,b]$ to $[b',c']$,
where $c>p_1$, and thus the new long zigzag ends in the next chamber
through which the embedding $h$ traverses. For $k>1$ we apply inductively the same procedure.   
\end{proof}
\begin{proof}[Proof of Lemma \ref{l:ample-parametric} for $q\geq 3$]
We begin with a set of fibered long  zigzags $\wt Z_j$ subordinated to 
a covering $ K=\mathop{\bigcup}\limits_1^N U_j$.  Given a point $s\in U_j$
for each $j=1,\dots, N$ we denote by  $Z^s_j$ the zigzag over $s\in U_j$,
and by $h^s_j:[a,b]\to V$   its parameterization.   

\n   
Fix a $j=1,\dots N$ and denote  by $C_0$  one of the chambers to which
the zigzag $\wt Z_j$ is adjacent,
say $h^s_j(b)\in C_0$ for $s\in U_j$. Passing, if necessary, to a set
of zigzags subordinated to a finer covering of $K$  we can  extend the family
of embeddings $h^s_j$, $s\in U_j$,  to a family of embeddings
$[a,c]\to V$, $c>b,$  still denoted by $h_j^s$,  such that 
\begin{itemize}
\item for each $ s\in U_j$ the embedding $h^s_j$
is transversal to $S$ and  $H_j^s(c)\in C$;
\item the image $f^s_j([a,c])\subset W$ is independent of $s$ and disjoint
from images of all    other zigzags.
\end{itemize}
Using Lemma \ref{l:moving-long}
we can construct  a deformation of   
the fibered map $\wt f|_{U_j}$  which is supported in 
$\Op\wt h_j([a,c])$ which creates a fibered zigzag adjacent to the chamber $C$.
Note that  for  different  $i=1,\dots, k$, and  $j=1,\dots, N$
all these deformations are supported in non-intersecting neighborhoods, 
and hence can be done simultaneously.
\end{proof}

{\sc  Case  of open $W$  and $q\geq2$.}
Let us consider a function $\phi:W\to\bbR$ without critical points. 
Let $\F$ be a $1$-dimensional foliation by its gradient trajectories of $\phi$ 
for some Riemannian metric.
 
\bp
\label{l:open-ample}
{\bf(The case $K=$ point)}
Any  $f\in\fM(V,W,S)$ admits a zigzag adjacent to any of  its chambers
which projects to one of the leaves of the foliation $\F$.
\ep

\begin{proof}
Let us call a leaf   $L$ of $\F$ {\it regular} for $f$ if
the map $f|_S:S\to  W$ is transversal to it.
Sard's theorem implies that  a generic leaf of $\F$ is regular. 
Let $C$ be one of the chambers. Given a regular for $f$   leaf $L$
we denote by $C_L$  the union of those connected components of the closed
1-dimensional manifold $f^{-1}(L)\subset V$ which intersect  $C$.
We claim that for some  regular leaf $L$   there exists
a non-empty component $\Sigma$ of $C_L$ such   that $f|_S:S\to L$ has 
more than 2 critical points. Indeed, otherwise we could reconstruct
an involution $\alpha:V\to V$ such that $f\circ \alpha=f$, which would  imply
that  $f$ is taut. 
The intersection of $C$ with the circle $\Sigma$ consists of one or several 
arcs.
Let $A$ be one of these arcs. Its end points, $p_1$ and $p_2$, belong to
the fold $S$, and are critical points of $f|_\Sigma:\Sigma\to L$.
Let us assume  that $p_1$ is a local minimum and $p_2$ is a local maximum.
Recall that the leaf $L$ is oriented by the gradient vector field of
the function $\phi$. We orient the arc $A$ from $p_1$ to $p_2$ and orient
the circle  $\Sigma$  accordingly. Let $p_3 $ be the next critical point of 
$f|_\Sigma$.
Choose  a point  $q_1\in A$ close to $p_1$ and a point $q_2$ close to $p_3$
and after it in terms of the orientation.  If $f(p_3)>f(p_1)$ then
the arc $Z=[q_1,q_2]$ is a zigzag adjacent to $C$ beginning at $q_1$
and ending at $q_2$. If  $f(p_3)<f(p_1)$ then the same arc $Z=[q_1,q_2]$
is again is a zigzag adjacent to $C$, but  beginning at $q_2$ and ending at 
$q_1$.
\end{proof}
\begin{proof}[Proof of Lemma \ref{l:ample-parametric} for an open $W$ and  
$q\geq 2$]
Lemma \ref{l:open-ample} implies that over any point $s\in K$ the   map $f_s$ 
admits
a zigzag  adjacent to the chamber $C$. Moreover, we can assume that all these 
zigzags
lie over different leaves of $\F$.  These zigzags extend to a neighborhood $\Op 
s\in K$
as fibered zigzags over this neighborhood   which for each $s'\in \Op s$ project
to the same leaves of $\F$. Hence, we can choose a finite covering 
$K=\mathop{\bigcup}\limits_1^N U_j$ by neighborhoods $U_j$ over which 
there exist ample sets of special  fibered zigzags.  Clearly we can arrange
that zigzags over different $U_j$ project to different leaves of $\F$, and thus 
their images do not intersect.
\end{proof}
\section{Proof of the main theorem}
\subsection{Wrinkled $S$-immersions}
A map  $f:V\to W$ is called a {\it  wrinkled $S$-immersion} if it
has $S$ as its  fold singularity, and in the complement of $S$ it is
a wrinkled map, see below Section \ref{ss:wrinkles}. Let $C$ be one of
the chambers (i.e. a connected
component of $V\setminus S$). We denote by $\fM_w(V,W,S,C)$   
the space of wrinkled $S$-immersions which have all wrinkles in the chamber $C$.
We will call $C$ the {\it designated} chamber. 
The regularized differential construction (see Section \ref{ss:wrinkles}) 
provides a map
$d_R:\fM_w(V,W,S,C)\to\fm(V,W,S)$, and the Wrinkling Theorem
\ref{t:B} implies that 
\bp
\label{t:formal-wrinkled}
{\bf (From formal  to wrinkled $S$-immersions)}
The map $d_R$ is a (weak) homotopy equivalence.
\ep
\begin{proof} 
Let $C_1, \dots C_l$ be a sequence of all chambers in $V\setminus S$
such that $C_l=C$ and each $C_{i}$, $i<l$,  have a common wall with $C_j$, 
$j>i$. 
We apply Hirsch's $h$-principle (see \cite{[Hi]}) for equidimensional
immersions of open (i.e. non-closed!)
manifolds to $\ol C_1$ and then make a fold on $\p \ol C_1$.
Next,  we apply the relative version the same $h$-principle
to the pair $(C_2,\p C_1)$
and so on. On the last step we apply the Wrinkling theorem
\ref{ss:wrinkles} to the pair $(\ol C,\p \ol C)$. 
\end{proof}

A little bit stronger statement can be formulated in the language of fibered 
maps.
A fibered  over $K$ map is called a {\it fibered  wrinkled $S$-immersion} if it 
has
$S$ as its  fiberwise fold  singularity, and in the complement of
$K\times S\subset K\times V$ it is a fibered wrinkled map.
One can also talk about fibered over $K$  formal  $S$-immersions,
i.e. parameterized by $K$ families $F_s\in\fm(V,W,S)$. 
Then a regularized differential of a fibered $S$-immersion is
a fibered formal $S$-immersion. The following proposition is a slight 
improvement of the above homotopy equivalence claim and it also follows from 
\ref{t:B}.

\bp
\label{t:formal-wrinklfibered}
{\bf (From formal  to wrinkled $S$-immersions;  a fibered version)}
Given any fibered  over $K$ formal $S$-immersion $\wt F\in \fm^K(V,W,S)$, there 
exists a fibered
over $K$  wrinkled map  $\wt g\in\fM^K_w(V,W,S,C)$ whose regularized fibered 
differential
is homotopic to $\wt F$. Moreover, if over a closed $L\subset K$ we have $\wt 
F=d\wt f$, 
where $\wt f$ is a genuine fibered $S$-immersion, then the map $\wt g$ can be 
chosen
equal to $\wt f$ over $L$, and the homotopy can be made fixed over $L$.
\ep
Now we want to supply a wrinkled $S$-immersion by zigzags adjacent to
the designated chamber. 
\bp
\label{p:wrinkled-zigzags}
{\bf (Zigzags for fibered wrinkled $S$-immersions)}
Let $\wt f$ be a   fibered over $K$ wrinkled map which have all
wrinkles in the chamber $C$.  Suppose that  over  a closed
subset $L\subset K$ the fibered map $\wt f$ consists of genuine (i.e non-
wrinkled) $S$-immersions.
Let $\,\mathop{\bigcup}\limits_1^N U_j\supset L$ be a covering of $L$ by 
contractible and open in $K$ sets, and $\wt Z_1,\dots,\wt Z_N$ be  a set of 
adjacent to $C$ fibered zigzags subordinated to the covering 
$\mathop{\bigcup}\limits_1^N U_j$.  
Then there is a homotopy
$\wt f_t\in\fM^K_w(V,W,S,C)$,  $t\in[0,1]$,  $\wt f_0=\wt f$,   such that
\begin{itemize}
\item $\wt f_t$ is fixed over $L$ in a neighborhoods of the zigzags 
$\wt Z_1,\dots,\wt Z_N$;
\item  there exists a zigzag $\wt Z$,  fibered over a domain  $U$,
$K\setminus\mathop{\bigcup}\limits_1^N U_j\subset U\subset K\setminus \Op L$,
for $\wt f_1$, adjacent to $C$ and disjoint from $\wt Z_1,\dots,\wt Z_N$.
\end{itemize}
\ep
Let us denote $\Sigma:=\p C\subset S$.  
The following two lemmas  will be needed in the proof of \ref{p:wrinkled-zigzags}. 
\bp
\label{l:disjoint-section}
Let  $U_1,\dots U_N$ be as in \ref{p:wrinkled-zigzags}. 
Let $\sigma_j:U_j\to U_j\times \Sigma$, $j=1.\dots, N$, be disjoint 
sections.
Then there exists a section $\sigma:K\to K\times \Sigma$
disjoint from the  sections $\sigma_1,\dots,\sigma_N$ and homotopic to a 
constant section.
\ep

\begin{proof} 
Arguing by induction over $j=1,\dots, N$, we construct a fiberwise isotopy $\wt 
g_t:K\times V\to V$, $t\in[0,1]$,  which makes sections $\sigma_j: U_j\to 
U_j\times \Sigma$ constant, i.e. $\wt g_1\circ\sigma_j(s)=(s,c_j)$, 
$c_j\in\Sigma$, $j=1,\dots, N$.
Take a point $c\in\Sigma$ different from $c_1,\dots, c_N$ and define 
$\sigma'(s):=(s,c)$ for any $s\in K$. Then the section
$\sigma=\wt g_1^{\,-1}\circ\sigma':K\to K\to V$ has the required properties.
\end{proof}

\bp 
{\bf(Local birth of a zigzag near $\sigma$)}
\label{p:local-zigzag}
Let  $M\subset K$ be a closed subset and 
$\s:M\to M\times \Sigma\subset K\times V$ a section.
Then there exists a deformation $\wt f_t\in\fM^K_w(V,W,S,C)$ of $\wt f_0=\wt f$ 
which is supported in 
$\Op \s(M)\subset K\times V$ such that $\wt f_1$ admits a fibered over $\Op M$ 
zigzag adjacent to $C$ and supported in $\Op \s(M)\subset K\times V$.
\ep

\begin{proof} 
First we give a sketch of the construction. 
Take a fibered embedding $\wt h:M\times I\to M\times\Op \Sigma$
such that for all $s\in M$ the image $h(s\times I)$
lies on a small $\alpha_{loc}$-invariant circle near $\sigma(s)\in\Sigma$,
see Fig.\ref{em28}. Then there exists a homotopy $\wt g_t$  of the map
$\wt g_0=\wt f|_{\Op h(M\times I)}$
such that $\wt g_1$
have a fibered zigzag over $M$,
see Fig.\ref{em28}. This local homotopy can be extended as a homotopy
$\wt f_t\in\fM^K_w(V,W,S,C)$ of the whole fibered map $\wt f$. 

\mn
Let us give now a more detailed description.
Let $S^1\subset \bbC$ be the unit circle, and $\exp:\bbR\to S^1$ 
be a covering map $u\mapsto e^{ iu}$, $u\in\bbR$.  Choose a neighborhood 
$\Omega\supset \s(M)$.
Let $\wt h:M\times S^1\to \Omega\subset M\times  V$ be a fibered 
embedding such that for all $s\in M$ the image $h^s(S^1)$
is  a small $\alpha_{loc}$-invariant circle near $\sigma(s)\in\Sigma$,
see Fig.\ref{em28}.     
\begin{figure}
\centerline{\psfig{figure=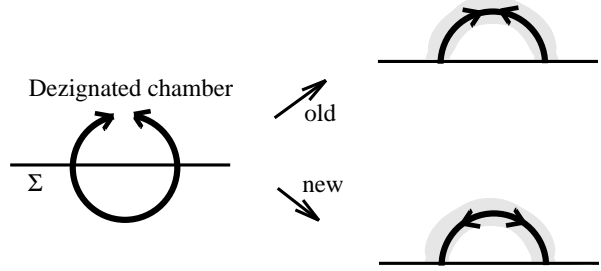,height=45mm}}
\caption{Local birth of a zigzag}
\label{em28}
\end{figure}
We can assume that $h^s(\exp u)\in \oC$ for $u\in[\pi/2,3\pi/2]$ and $h^s(\exp 
u) \notin \oC$ for $u\in(-\pi/2,\pi/2)$.
Consider also  a fibered embedding 
$$\wt\varphi:\Op M\times B\to \Op M\times C\subset \Omega\,, $$
where $B$ is an open 
$n$-ball, such that  $\varphi^s(B)$ is a small ball centered at the point $h^s(-
1)\in C$ and such that 
$\varphi^s(D)\cap h^s(S^1)=h^s(\exp((\pi-\e,\pi+\e))$,   $s\in M$.
There exists a  compactly supported fibered regular  homotopy,  see 
Fig.\ref{em28},
$$\wt \psi_t:\Omega\setminus \wt \varphi(\Op M\times B)\to \Omega,\;\;t\in[0,1],$$ 
such that $\wt\psi_t(\wt h(s,\exp u))=\wt h(s,\exp(1+\frac {t}2)u)$,
$s\in M$, $u\in[-\pi+\e,\pi-\e]$. Notice that  $ \psi^s_1$ maps an embedded  arc 
$h^s([-\pi+\e,\pi-\e])$ onto an overlapping arc $h^s([-\frac32(-\pi+\e), 
\frac32(\pi-\e))$. 
The regular homotopy  $\wt\psi_t$ extends, according to Theorem  \ref{t:B}, 
to a compactly fibered wrinkled homotopy $\Omega\to\Omega$, 
and then can be extended further to  the rest of $K\times V$ as the identity 
map. We will 
use the same notation $\wt\psi_t$ for this extension.  Finally,    
the   wrinkled homotopy $\wt f_t=\wt f_0\circ\wt\psi_t:K\times V\to K\times V$ 
connects $\wt f_0$ with a wrinkled map $\wt f_1$  which has the embedding 
$\wt\psi_1\circ\wt h|_{M\times [-\pi-\e,\pi+\e])}$ as its fibered over $M$ 
zigzag.
\end{proof}
\begin{proof}[Proof of Proposition \ref{p:wrinkled-zigzags}]
Note that the zigzags
$\wt Z_j$ over  $U_j\subset K$, $j=1,\dots, N$, intersect the closure
$\oC$ of the chamber $C$ along intervals with one  end on $\Sigma$ and
the second inside $C$.   Taking the end-points in $\Sigma$ we get  sections 
$\sigma_j:U_j\to U_j\times\Sigma$,
$j=1,\dots, N$. Let us  apply Lemma \ref{l:disjoint-section} and construct a 
section $\sigma:K\to K\times\Sigma$ disjoint from the sections 
$\sigma_1,\dots,\sigma_N$.
Let $M\subset K\setminus L$ be a  closed set
such that $K\setminus M\subset \mathop{\bigcup}\limits_1^N U_j$. There exists a 
neighborhood   $\Omega\supset \s(M)\subset K\times V$ which does not intersect 
the fibered zigzags $\wt Z_j$, $j=1,\dots, N$. Let $U\supset M$ be an open 
neighborhood whose closure  is contained in $K\setminus L$, and such that 
$\sigma(U)\subset \Omega$.    We  conclude the proof of
\ref{p:wrinkled-zigzags} by applying
Lemma \ref{p:local-zigzag}.
\end{proof}

\subsection{Engulfing wrinkles by zigzags}
\label{ss:engulfing}
\bp
\label{t:engulfing}
{\bf (Getting rid of wrinkles)}
Let $\wt f\in\fM^K_w(V,W,S,C)$ admits a set of fibered zigzags adjacent
to the chamber $C$. Suppose that $\wt f$   is a genuine fibered $S$-immersion over a 
closed $L\subset 
K$.
Then $\wt f$ is homotopic in $\fM^K_w(V,W,S,C)$ to a genuine fibered $S$-
immersion via a homotopy fixed over $L$.
\ep
\begin{proof}
Let   $\wt Z_1,\dots,\wt Z_N$ be fibered zigzags adjacent to $C$ and 
subordinated
to a covering  $\mathop{\bigcup}\limits_1^N U_j=K$.
First of all, we can apply the enhanced wrinkling theorem (see the remark  after 
Theorem \ref{t:B})
to $\wt f|_{K\times C}$ and get a modified $\wt f$ such that each fibered 
wrinkle 
is supported over one of the elements of the covering. Using Lemma 
\ref{l:inscribed} we can assume that wrinkles and  zigzags are
in  a 1-1 correspondence with the elements of the covering.
We will eliminate inductively all the wrinkles by a procedure which we call
{\it engulfing of wrinkles by zigzags}.
We will  discuss   this construction only in the non-parametric case,  i.e. when 
$K$ is a  point. The case of a general $K$ differs only in the notation.

\mn
Let $w=w(q)$ be the standard wrinkle with the membrane
$D^q$ (see \ref{ss:wrinkles}).
Let us recall that $w$ is  fibered over $\bbR^{q-1}$ 
map $\bbR^q\to\bbR^q$ defined by the formula
$$
(y,z) \mapsto\left( y,z^3 + 3(|y|^2-1)z \right),
$$
where $y\in \bbR ^{q-1},\,z\in \bbR^1$, and
$|y|^2 = \displaystyle{\sum_1^{q-1}y^2_{i}}$. 
Note that  for a sufficiently small $\e>0$ we have  
$w(D^q_{1+\e})\subset \{|y|\leq 1+\e, |z| \in[-3,3]\}=D^{q-1}_{1+\e}\times[-3,3]$.

We will need the following lemma.

\bp
\label{l:st-engulf}
{\bf (Standard model for engulfing)}
Let us consider a fibered over $D^{q-1}_{1+\e}$,  map $\Gamma: D_{1+\e}^{q-
1}\times[a,c] \to D_{1+\e}^{q-1} \times [a,c]$   such that
\begin{itemize}
\item  for some $b\in(a,c)$  the restriction $$\Gamma|_{D^{q-
1}_{1+\e}\times[a,b]\ }: D^{q-1}_{1+\e}\times[a,b]\to  D^{q-
1}_{1+\e}\times[a,c]$$ is a fibered over $D^{q-1}_{1+\e}$ zigzag;
\item   the restriction 
$$\Gamma|_{D^{q-1}_{1+\e}\times[b,c]}:D^{q-1}_{1+\e}\times[b,c]\to D^{q-1}_{1+\e}\times[a,c]$$ 
is a fibered over $D^{q-1}_{1+\e}$ wrinkled map with the unique wrinkle 
whose cusp locus projects to the sphere
$\p D^{q-1}$. In other words, there exist a fibered over $D^{q-1}_{1+\e}$  
embeddings  $\alpha:D^q_{1+\e}\to D^{q-1}_{1+\e}\times[b,c]$  and 
$\beta:D^{q-1}_{1+\e}\times[-3,3]\to D^{q-1}_{1+\e}\times[a,c] $ such that 
$\Gamma\circ\alpha=\beta\circ w$.
\end{itemize}
Then there exists a fibered over $D^{q-1}_{1+\e}$ homotopy 
$$\Gamma_t: D_{1+\e}^{q-1}\times[a,c] \to D_{1+\e}^{q-1} \times [a,c],\, t\in[0,1]\,,$$   
which begins with $\Gamma_0=\Gamma$,  is fixed near
$D_{1+\e}^{q-1}\times[a,c]$ and  such that
\begin{itemize} 
\item $\Gamma_t|_{D_{1+\e}^{q-1}\times[a,b] }$ is the homotopy in the space of 
fibered folded maps;
\item   $\Gamma_t|_{D_{1+\e}^{q-1}\times[b,c] }$ is a fibered wrinkled 
homotopy which eliminates the wrinkle of the map $\Gamma_0=\Gamma$, i.e. the map  
$\Gamma_1|_{D_{1+\e}^{q-1}\times[b,c] }$  is non-singular.
\end{itemize}
\ep
We will call the homotopy $\Gamma_t$ {\it engulfing of the wrinkle $w$ by a zigzag.}
\begin{proof} 
For $y\in \Int\, D^{q-1}$ the homotopy $\Gamma_t$ is shown on
Fig.\ref{em23}. This deformation can be done smoothly depending on the parameter 
$y\in D^{q-1}_{1+\e}$ and dying out  on $[1,1+\e]$. See also Fig.\ref{s24}, 
where the deformation $\Gamma_t$ is shown for $q=2$.
\end{proof}
 
\begin{figure}
\centerline{\psfig{figure=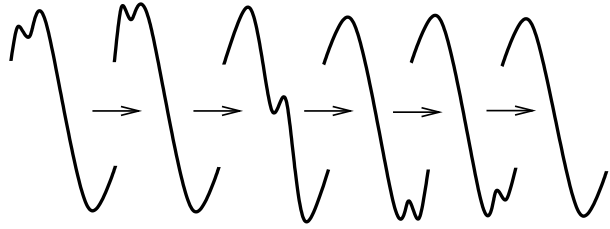,height=30mm}}
\caption{Homotopy ${\Gamma_t}$ on $y\times [a,c]$, $y\in \Int\,D^q$}
\label{em23}
\end{figure}

\begin{figure}
\centerline{\psfig{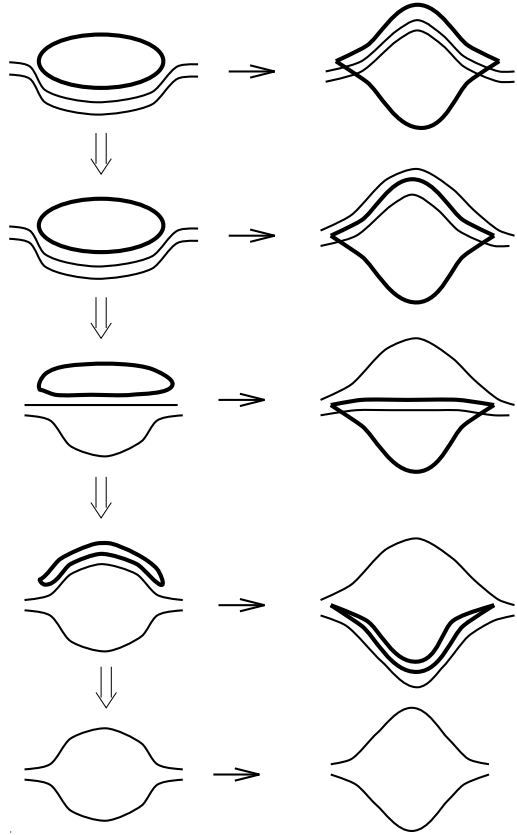}}
\caption{Homotopy $\Gamma_t$, $q=2$}
\label{s24}
\end{figure}

We continue now the proof of \ref{t:engulfing} by reducing it to the standard 
engulfing model \ref{l:st-engulf}.
Let us recall that there is 1-1-correspondence between 
the wrinkles and zigzags.
Let $h:[a,b]\to V$ be a zigzag of a map $f$ adjacent to the designated  chamber  
$C$. The embedding $h$ extends to an embedding
$H:D_{1+\e}^{q-1}\times[a,b]\to V$  onto a small neighborhood of
the zigzag $Z=h([a,b])$, such that
$H|_{y\times [a,b] }$ is a zigzag for each $y\in D^{q-1}$.
Take the  corresponding wrinkle with a membrane $D\subset C$. By definition 
this means that for a sufficiently small $\e>0$ there exists an embedding 
$\alpha: D^q_{1+\e}\to V$ such that $\alpha (D^q)=D$ and $f\circ\alpha 
=g\circ w(q)$ for an embedding $g:D^{q-1}_{1+\e}\times[-3,3]\to W$.  

As it  is clear from Fig.\ref{em271}, there exists  an extension of the  
embedding  $H$ to an embedding $\wh H: D_{1+\e}^{q-1}\times[a,c]\to V$, $c>b$,
such that
\begin{itemize}
\item   $\wh H(D_{1+\e}^{q-1}\times[b,c])\supset D$;
\item  the map $f\circ \wh H: D_{1+\e}^{q-1}\times [a,c]$ can be written as  
$g\circ \Gamma$, where $\Gamma$ is a fibered over $D^{q-1}_{1+\e}$ map $ 
D_{1+\e}^{q-1}\times [a,c]\to D_{1+\e}^{q-1}\times [a,c]$ and $g: 
D_{1+\e}^{q-1}\times [a,c]\to W$ is an immersion;
\item  the restriction of $\Gamma$ to   $D_{1+\e}^{q-1}\times [b,c]$ is a 
fibered over $D^{q-1}_{1+\e}$ wrinkled map with the unique wrinkle whose    
cusp locus  projects to the sphere $\p D^{q-1}\subset D^{q-1}_{1+\e}$.
\end{itemize}
\begin{figure}
\centerline{\psfig{figure=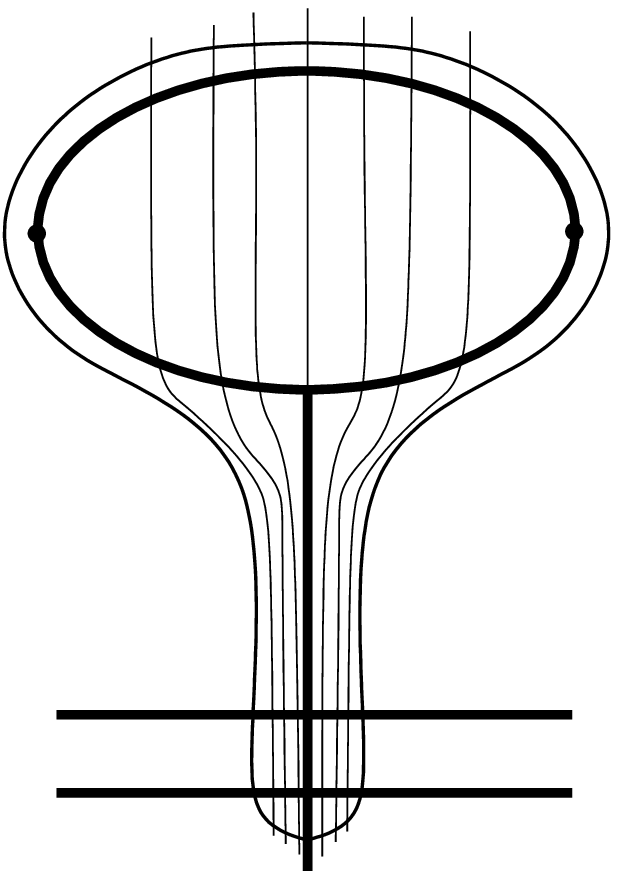,height=60mm}}
\caption{Embedding  $\wh H$}
\label{em271}
\end{figure}

Thus we  are in a position to apply Lemma \ref{l:st-engulf}.   Let  
$$\Gamma_t:[a,c]\times D_{1+\e}^{q-1}\to [a,c]\times D_{1+\e}^{q-1},\;t\in[0,1],$$ 
be the engulfing homotopy  of $\Gamma=\Gamma_0$ constructed  in
\ref{l:st-engulf}, which eliminates the wrinkle.   Note that  the homotopy $\Gamma_t$ is 
fixed near $\p\left([a,c]\times D_{1+\e}^{q-1}\right)$. 
This enables us to define the required wrinkled  homotopy  
$f_t:V\to W$, $t\in[0,1]$, which eliminates the wrinkle $w$ by setting it 
equal  to   $f$ on $V\setminus\wh H([a,c]\times D^{q-1}_{1+\e})$, and to 
$g\circ\Gamma_t\circ\wh H^{-1}$ on $\wh H([a,c]\times D^{q-1}_{1+\e})$.
\end{proof}

\subsection{Proof of Theorem \ref{t:soft}}
 Let us prove that the map $d$ induces an injective homomorphism  on $\pi_{k-
1}$, i.e.
  that $\pi_k({\frak m}(V,W,S),\,{\frak  M}_{soft}(V,W,S))=0$.
Denote $K:=D^k$, $L:=\p D^k=S^{k-1}$.  
We need to show that if $\wt F\in\fm^K(V,W,S)$ is a fibered formal $S$-
immersion, 
which is equal to $d\wt f$ over $L$ for $\wt f\in \fM^L_{soft}(V,W,S)$, 
then there exists a homotopy $\wt F_t$ of  $\wt F$ in 
$(\fm^K(V,W,S),\fM^L_{soft}(V,W,S)$
such that $\wt F_1=d\wt f_1$ for  $\wt f_1\in  \fM^K_{soft}(V,W,S)$.
The construction of the required homotopy can be done in four steps: 
\bi
\item Using Lemma 
\ref{l:ample-parametric}  we first  deform $\wt F$ in $(\fm^K(V,W,S), 
\fM^L(V,W,S))$
so  that the new $\wt F$ admits over $L$ a set of special fibered zigzags
adjacent to a designated chamber $C$.
\item Using  \ref{t:formal-wrinkled} we   further deform $\wt F$ 
into a differential of a fibered wrinkled $S$-immersion  $\wt f$ such that all 
the  
wrinkles belong to the chamber $C$. 
\item Using \ref{p:wrinkled-zigzags}  we  can further deform $\wt f$  in such a 
way  that 
the new $\wt f$ admits a set of zigzags adjacent to the chamber $C$. 
Furthermore,  we can choose this set of zigzag in such a way that it includes 
the 
set of zigzags over $L$.
\item Using \ref{t:engulfing} one can engulf all the fibered wrinkles.
\ei
 The proof of the surjectivity claim is similar but does not use the first step.

\section {Appendix: folds, cusps and wrinkles}
\label{s:wrinkles}

We recall here, for a convenience of the reader,
some definitions and results from \cite{[EM97]}. 
We consider here only the equidimensional case $n=q$, 
and this allows us to  simplify the definition, the  notation etc.,  compared to 
\cite{[EM97]}.     

\subsection{Folds and cusps}
\label{ss:folds}

Let $V$ and $W$ be smooth manifolds of the same  dimension  $q$.
For a smooth map $f:V \rightarrow W$ we will denote by
$\Sigma(f)$ the set of its singular points, i.e.
$$
\Sigma(f) = \left\{ p \in V, \; \hbox { rank } d_pf < q \right\}\;.
$$
A point $p \in \Sigma(f)$ is called a {\it fold} type singularity or
a {\it fold} of index $s$ if near the point $p$ the map $f$ is equivalent to
the map
$
\bbR ^{q-1} \times \bbR ^{1} \rightarrow \bbR ^{q-1}
\times \bbR ^{1}
$
given by the formula
$
(y,x) \to ( y,\,\, x^2)
$
where $y=(y_1,...,y_{q-1})\in \bbR ^{q-1}$.

Let $q>1$. A point $p \in \Sigma(f)$ is called a {\it cusp}
of index $s+\frac{1 }{ 2}$\; if near the point $p$ the map $f$ is
equivalent to the map
$
\bbR ^{q-1}\times \bbR ^1\rightarrow \bbR^{q -1}\times \bbR ^1
$
given by the formula
$
(y,z) \to ( y,z^3 + 3y_1z )
$
where $z \in \bbR ^1,\,\,
y = (y_1, \ldots, y_{q-1}) \in \bbR ^{q-1}$.

\mn
For $q\geq 1$ a point $p \in \Sigma(f)$ is called an 
{\it embryo} of index $s + \frac{1}{ 2}$\,\, if $f$
is equivalent near $p$ to the map
$
\bbR^{q-1}\times\bbR^1\rightarrow
\bbR^{q-1}\times\bbR ^1
$
given by the formula
$
(y,z) \to (y, z^3 + 3|y|^2z)
$
where $y \in \bbR ^{q-1}$,\,\, $z \in \bbR ^1,\,\,
|y|^2=\sum\limits_1^{q-1}y_i^2$.
The set of all folds of $f$ is denoted by $\Sigma^{10}(f)$, the set of 
cusps by
$\Sigma^{11}(f)$ and the closure $\overline{\Sigma^{10}(f)}$ by
$\Sigma^1(f)$.

\mn
Notice that folds and cusps are stable singularities for individual
maps, while embryos are stable singularities only for $1$-parametric
families of mappings. For a generic perturbation of an individual map
embryos either disappear or give birth to wrinkles which we consider
in the next section.

\subsection{Wrinkles and wrinkled maps}
\label{ss:wrinkles}

Consider the map
$
w(q): \bbR^{q-1}  \times \bbR^1
\rightarrow \bbR^{q-1} \times \bbR^1
$
given by the formula
$$
(y,z) \mapsto
\left( y,z^3 + 3(|y|^2-1)z\right),
$$
where $y \in \bbR ^{q-1}, \,z \in \bbR^1$ and
$|y|^2 = \displaystyle{\sum_1^{q-1}y^2_{i}}$.

\begin{figure}
\centerline{\psfig{figure=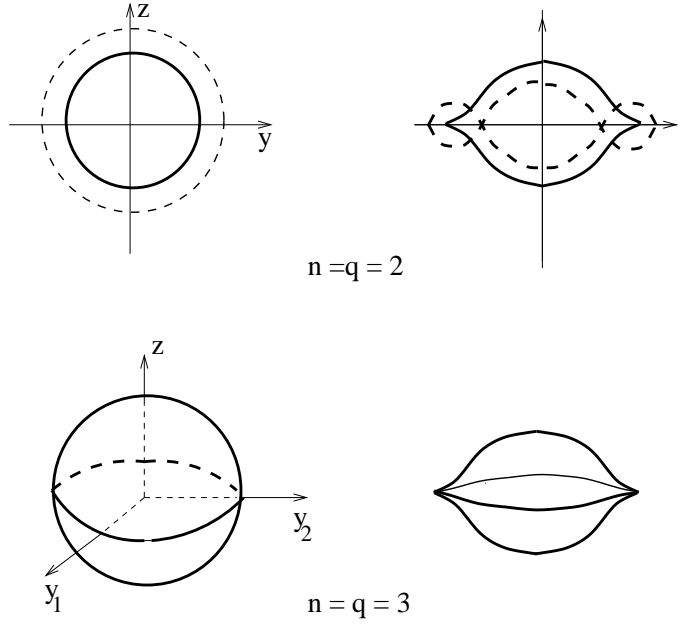,height=85mm}}
\caption{Wrinkles. Pictures in the source and in the image.}
\label{em5}
\end{figure}

\n
The singularity
$\Sigma^1(w(q))$ is the $(q-1)$-dimensional sphere
$
S^{q-1}\subset\bbR^q.
$
Its equator $\{z=0, |y|=1 \} \subset \Sigma^1(w(q))$
consists of cusp
points, the upper and lower hemisphere
consists of folds points (see Fig.\ref{em5}). We will call the $q$-dimensional
bounded by $\Sigma^1(w)$
disc $D^q=\{z^2 + |y|^2\leq 1\}$ 
the {\it membrane} of the wrinkle.

\mn
A map $f: U \rightarrow W$ defined on an open subset $U\subset V$
is called a {\it wrinkle} if it is equivalent to the
restriction of $w(q)$ to an open neighborhood $U$
of the disk $D^q\subset \bbR^q$.
Sometimes  the term
``wrinkle''  can be also used also for the singularity $\Sigma (f)$ of  the wrinkle $f$.

\mn
Notice that for $q=1$ the wrinkle is a function with two nondegenerate
critical points of indices $0$ and $1$ given in a neighborhood of a gradient
trajectory which connects the two points.

\mn
Restrictions of the map $w(q)$ to subspaces $y_1=t$,
viewed as maps $\bbR ^{q-1} \to \bbR ^{q-1}$, are non-singular maps
for $|t|>1$, equivalent to $w(q-1)$ for $|t|<1$ and to embryos for
$t=\pm 1$.

\mn
Although the differential $dw(q):T(\bbR^q)\rightarrow T(\bbR^q)$
degenerates at points of $\Sigma(w)$, it can be canonically {\it regularized}.
Namely, we can change the element $3(z^2 + |y|^2 -1)$ in the Jacobi matrix of
$w(q)$ to a function $\gamma$ which
coincides with $3(z^2+|y|^2-1)$ outside an arbitrary small
neighborhood $U$ of the disc $D^q$ and does not vanish on $U$.
The new bundle map
${\cal R}(dw):T (\bbR^q)\rightarrow T(\bbR^q)$
provides a homotopically canonical
extension of the map $dw:T(\bbR ^q\setminus U)\rightarrow T(\bbR ^q)$ to
an epimorphism (fiberwise surjective bundle map)
$T(\bbR^q)\rightarrow T(\bbR^q)$. We call
${\cal R}(dw)$ the {\it regularized differential} of the map $w(q)$.

\mn
A map $f:V \rightarrow W$ is called {\it wrinkled} if there exist disjoint
open subsets $U_1, \ldots, U_l \subset V$ such that $f|_{V\setminus  U}, \,\,
U=\bigcup^{\,l}_1 U_i,$ is an
immersion (i.e. has rank equal $q$) and for each $i=1,\ldots,l$ the 
restriction $f|_{U_i}$ is a wrinkle. Notice that the sets $U_i,\,i=1,\dots,l,$
are included into the structure of a wrinkled map.

\mn
The singular locus $\Sigma\,(f)$ of a wrinkled map $f$ is a union of
$(q-1)$-dimensional wrinkles $S_i=\Sigma^1(f|_{U_i})
\subset U_i$.
Each $S_i$ has a $(q-2)$-dimensional equator $T_i \subset S_i$
of cusps which divides $S_i$ into 2 hemispheres of folds of 2 neighboring
indices. The differential $df:T(V)\rightarrow T(W)$ can be regularized to
obtain an epimorphism
$\R (df):T(V) \rightarrow T(W)$. To get
$\R (df)$ we regularize $df|_{U_i}$ for each wrinkle $f|_{U_i}$.

\subsection{Fibered wrinkles}\label{ss:fibrinkles}

All the notions from \ref{ss:wrinkles} can be extended
to the parametric case.

\mn
A {\it fibered} (over $B$) {\it map} is a commutative diagram
$$
\xymatrix{X \ar[rd]_p\ar[rr]^f& &Y\ar[ld]^q \\&B&}
$$
where $p$ and $q$ are submersions. A fibered map can be also denoted simply by 
  $f:X\rightarrow Y$ if $B, p$ and $q$ are implied from the context.

\mn
For a fibered map
$$
\xymatrix{X \ar[rd]_p\ar[rr]^f& &Y\ar[ld]^q \\&B&}
$$
we denote by $T_BX$ and $T_BY$ the subbundles Ker $p \subset TX$ and Ker $q
\subset TY$.  They are tangent to foliations of $X$ and $Y$ formed by
preimages
$$p^{-1}(b) \subset X\,,\,\, q^{-1}(b) \subset Y\,,\,\, b \in B\,.$$

\mn
The fibered homotopies, fibered differentials, fibered submersions, and so on
are naturally defined in the category of fibered maps (see \cite{[EM97]}).
For example, the {\it fibered differential} of
$f:X \rightarrow Y$ is the restriction
$$d_Bf = df|_{T_BX} : T_{B}X \rightarrow T_BY\,.$$
Notice that $d_Bf$ itself is  a map fibered over $B$.

\mn
Two fibered maps, $f:X\to Y$ over $B$  and $f:X'\to Y'$ over $B'$,
are called {\it equivalent} if there exist open subsets
$A \subset B,\,\,
A'\subset B',\,\,
Z\subset Y,\,\,
Z'\subset Y'$
with
$f(X)\subset Z,\,\,
p(X)\subset A,\,\,
f'(X')\subset Z',\,\,
p'(X') \subset A'$
and diffeomorphisms
$\varphi: X \rightarrow X', \; \psi : Z \rightarrow Z',\; s:A\rightarrow A'$
such that they form the following commutative diagram
$$
\xymatrix@C=14pt{X\ar[rrrrdddd]_p\ar[rrrd]^\varphi\ar[rrrrrrrr]^f&&&&&&&&Z
\ar[lllldddd]^q\ar[llld]_\psi\\&&&X'\ar[rd]_{p'}\ar[rr]^{f'}& &Z'
\ar[ld]^{q'}&&&\\&&&&A'&&&&\\&&&&&&&&\\&&&&A\ar[uu]_s&&&&}
$$
For any integer $k>0$ the map $w(k+q)$
can be considered as a fibered map $w_k(k+q)$ over
$\bbR^k\times {\bf 0}\subset \bbR^{k+q}$.
A fibered map equivalent to the restriction of $w_k(k+q)$
to an open neighborhood $U^{k+q}\supset D^{k+q}$ is called
a {\it fibered wrinkle}.
The regularized differential ${\cal R}(dw_k(k+q))$ is a
fibered (over $\bbR^k$) epimorphism
$$
\xymatrix
{\bbR^k \times T(\bbR^{q-1}\times\bbR^1)
\ar[rrr]^{\,\,\,\,\,\,\,\,{\cal R}(dw_k(q))}&&&\bbR^k
\times T(\bbR^{q-1} \times \bbR^1)}
$$
A fibered map $f:V \to W$ is called a {\it fibered wrinkled map}
if there exist disjoint open sets $U_1, U_2,\ldots,U_l\subset V$,
such that $f|_{V\setminus U},\,\, U=\bigcup^{\,l}_{\,1} U_i\,$
is a fibered submersion and for each $i=1,\ldots,l$ the
restriction $f|_{U_i}$ is a fibered wrinkle. The restrictions of
a fibered wrinkled map to a fiber may have, in addition to wrinkles,
embryos singularities.

\mn
Similarly to the non-parametric case one can define the
regularized differential of a fibered over $B$ wrinkled map $F:V\to W$,
which is a fibered epimorphism $\R(d_Bf):T_BV\to T_BW$.

\subsection {The Wrinkling Theorem}
\label{ss:main-wrinkling}
The following  Theorem \ref{t:A}
and its parametric version \ref{t:B} is the adaptation of the 
results of our paper \cite{[EM97]}  to  the simplest case $n=q$.
In fact in this case the results below can be also deduced from
a theorem of V. Po\'enaru, see \cite{[Po]}.

\bp 
{\bf (Wrinkled mappings)}
\label{t:A}
Let $F: T(V) \rightarrow T(W)$ be an
epimorphism which covers a map $f: V \to W$. 
Suppose that $f$ is an immersion on a neighborhood of a
closed subset $K \subset V$, and $F$
coincides with $df$ over that neighborhood.
Then there exists a wrinkled map $g: V\to W$ which coincides with
$f$ near $K$ and such that ${\cal R}(dg)$ and $F$ are homotopic
rel. $T(V)|_{K}$. Moreover, the map $g$ can be chosen arbitrary
$C^0$-close to $f$ and with arbitrary small wrinkles.
\ep

\bp
{\bf (Fibered wrinkled mappings)}
\label{t:B}
Let $f: V\to W$ be a fibered over $B$ map covered by a fibered
epimorphism $F: T_B(V) \to T_B(W)$.
Suppose that $f$ is a fibered immersion on a neighborhood of a
closed subset $K\subset V$, 
and $F$ coincides with $df$ near a closed subset $K \subset V$.
Then there exists a fibered wrinkled map
$g: V\to W$ which extends $f$ from a neighborhood of $K$, and such that the
fibered epimorphisms ${\cal R}(dg)$ and $F$ are homotopic 
rel. $T_B(V)|_K$. Moreover, the map $g$ can be chosen arbitrary
$C^0$-close to $f$ and with arbitrary small fibered wrinkles.
\ep

\n
{\bf Remark.} The proof (see \cite{[EM97]}) gives also
the following useful enhancement for \ref{t:B}:
{\it given an open covering $\{U_i\}_{i=1,\dots,N}$
of $\,B$, one can always choose $g$ such that for each fibered wrinkle 
there exists $U_i$ such that $p(D^{k+q})\subset U_i$,
where $D^{k+q}$ is the membrane of the wrinkle.}

\end{document}